\documentclass[12pt]{article}

\usepackage{amsmath,amsthm,amssymb,amsopn,amsfonts,pdfpages,url}
\usepackage[linesnumbered,ruled,vlined]{algorithm2e}
\usepackage{graphics,relsize}
\usepackage{graphicx}
\usepackage{multirow}
\usepackage{subfig}
\usepackage{bbm}
\usepackage{authblk}

\newcommand{\twospace}{\renewcommand{\baselinestretch}{1.3}\normalsize}

\newcommand{\Pry}{\mathbb{P}}

\newcommand{\e}{\mathbb{E}}
\newcommand{\R}{\mathbb{R}}

\newcommand{\rT}{\varrho_T}
\newcommand{\st}{\mathfrak{str}}
\newcommand{\dens}{\mathfrak{d}}
\newcommand{\denscu}{\mathfrak{d}_{X \cup Y}}
\newcommand{\densca}{\mathfrak{d}_{X \cap Y}}
\newcommand{\para}{{\mathcal R}}
\newcommand{\paraoo}{{\mathcal R}^{o}}
\newcommand{\x}{\mathcal{X}}
\newcommand{\HH}{\mathcal{H}}
\newcommand{\vecf}{\vec{\textup{f}}}
\newcommand{\vx}{X}
\newcommand{\vy}{Y}
\newcommand{\stuff}{\{0,\star,1\}^N}
\newcommand{\bS}{\overline{S}}
\newcommand{\bT}{\overline{T}}

\newtheorem{theorem}{Theorem}

\newtheorem{lemma}[theorem]{Lemma}
\newtheorem{corollary}[theorem]{Corollary}
\newtheorem{remark}[theorem]{Remark}
\newtheorem{conjecture}[theorem]{Conjecture}

\theoremstyle{definition}
\newtheorem{example}{Example}
\newtheorem{definition}{Definition}
\twospace
\begin{document}
\title{On a complete and sufficient statistic for the \\ correlated Bernoulli random graph model}
\author{Donniell E. Fishkind, Avanti Athreya, Lingyao Meng, \hspace{2in} Vince Lyzinski, Carey E. Priebe}
\maketitle
\begin{abstract}
Inference on vertex-aligned graphs is of wide theoretical and practical importance.
There are, however, few flexible and tractable statistical models for correlated
graphs, and even fewer comprehensive approaches to parametric inference on data
arising from such graphs. In this paper, we consider the correlated Bernoulli
random graph model (allowing different Bernoulli coefficients and edge correlations
for different pairs of vertices), and we introduce a new variance-reducing
technique---called \emph{balancing}---that can refine estimators for model parameters.
Specifically, we construct a disagreement statistic and show that it
is complete and sufficient; balancing can be interpreted as
Rao-Blackwellization with this disagreement statistic.
We show that for unbiased estimators of functions of model parameters,
balancing generates uniformly minimum variance unbiased estimators (UMVUEs).
However, even when unbiased estimators for model parameters do {\em not}
exist---which, as we prove, is the case with both the heterogeneity correlation
and the total correlation parameters---balancing is still useful, and lowers mean squared error.
In particular, we demonstrate how balancing can improve the efficiency
of the alignment strength estimator for the total correlation, a parameter that
plays a critical role in graph matchability and graph matching runtime complexity.\\

\noindent {\bf Keywords:} Rao-Blackwell, Lehmann-Scheff\' e, total correlation (graph), alignment strength, graph matching.

\end{abstract}

\newpage

\section{Overview \label{S:overview} } %%%%%%%%%%%%%%%%%%%%%%%%%%%%%%%%%

Paired random graphs with a natural alignment between their vertex sets arise in a
wide variety of application domains; for example, the interaction dynamics of the same
set of users across two social media platforms, or a pair of connectomes (brain graphs)
as imaged from two different subjects of the same species.
Given a pair of such graphs, the problem of {\em graph matching}---that is,
optimally aligning the two vertex sets in order to minimize edge disagreements, usually
with the purpose of obtaining the natural alignment---has a rich mathematical history,
and graph matching now plays a fundamental role in algorithms for machine learning and
pattern recognition; see the excellent surveys in \cite{thirtyyears, tenyears, gmsurvey}.

The {\it correlated Bernoulli random graph model}, described in Section \ref{sec:model},
is the focus of our work in this paper.
It is a versatile model used to describe two graphs that are correlated with each other
across a natural alignment between their vertex sets. The model allows for different
probabilities of adjacency for different pairs of vertices, and allows for different
edge correlations between different pairs of vertices across the natural alignment. This model is simple
enough to be theoretically and computationally tractable, yet it is rich enough to
successfully describe real data, and it has been profitably employed in this capacity; see, for example, \cite{SGM1,GMrelax,SGM2}.

{\bf The contributions of this paper fall into two groups, the second group utilizing the machinery
of the first group.}

{\bf Our first group of contributions:} In the context of a correlated Bernoulli random graph model,
we introduce a ``smoothing" procedure, called \emph{balancing}, which reduces the
mean-squared error for any estimator of a function of model parameters; specifically,
for any estimator $S$ of a function of model parameters $g(\theta)$, the balanced
estimator $\bS$ has the same bias as $S$, but has lower variance.
Indeed, under a  nondegeneracy condition, we prove in Theorem \ref{thm:second} that if $S$ is
an unbiased estimator of $g(\theta)$ then $\bS$ is the UMVU estimator of $g(\theta)$; this is
because $\bS$ is a Rao-Blackwellization of $S$ via the {\it disagreement statistic} $\HH$, and
we prove in our main result Theorem \ref{lem:complete} that $\HH$ is complete and sufficient, under the nondegeneracy condition.
We also prove in Theorem~\ref{thm:first}, under the  nondegeneracy condition, that if $S$ is
an unbiased estimator of $g(\theta)$ then any statistic
$T$ is also an unbiased estimator of $g(\theta)$ if and only if $\bS=\bT$.

These results should not be taken for granted; in Example \ref{ex:counter} we
illustrate that even knowing, hence fixing, the mean of the adjacency probabilities creates
a violation of the nondegeneracy condition, and the conclusions of the above theorems
then will indeed fail, in general.

Our second group of contributions of this paper focuses on
very recent advances in \cite{TotCor} regarding the correlated Bernoulli random graph model.
Specifically, the paper \cite{TotCor} introduced and showed the importance of a novel
model parameter called {\it total correlation},
which combines inter- and intra-graph contributions to a unified measure
of the correlation between the pair of graphs.
The authors empirically demonstrated---in broad families within the model---that graph matching complexity
and matchability are each functions of total correlation. They also proved that
the statistic called {\it alignment strength} is a strongly consistent estimator of total correlation.

{\bf Our second group of contributions:} In the context of a correlated Bernoulli random graph model,
the alignment strength statistic $\st$ was shown in \cite{TotCor} to be a strongly consistent
estimator of total correlation $\rT$ between the pair of graphs; however, we point out here that $\st$ is
{\bf not a balanced statistic}, hence, as noted above, the mean squared error in estimating $\rT$
is reduced by using $\overline{\st}$ instead.
We then prove (in Theorems \ref{thm:hom}, \ref{thm:tot}, \ref{thm:edg}) that there do not
exist unbiased estimators for several correlation parameters, including $\rT$.
Empirical experiments in Section~\ref{sec:empir} suggest that balancing the numerator and
denominator of $\st$ separately, which we call the {\it modified alignment strength} $\st'$, often has less
bias than $\overline{\st}$ as an estimator of $\rT$, always has less variance than $\overline{\st}$,
and we conjecture that $\st'$ always has
less mean square error~than~$\overline{\st}$~in~estimating~$\rT$.

The organization of this paper is as follows. The correlated Bernoulli random graph model,
important functions of the parameters, and important statistics are
described in Section \ref{sec:model}.
Our main results are stated in Section \ref{sec:results}
and proved in Section \ref{sec:pf}.
Empirical demonstrations~are~in~Section~\ref{sec:empir}.

\section{Correlated Bernoulli Random Graphs   \label{sec:model}}%%%%%%%%%%%%%%%%%%%%%
We begin by describing the correlated Bernoulli random graph model.
It consists of a pair of random graphs; without loss of generality these graphs are on the same
vertex set. (Indeed, the natural alignment between their vertex sets is a bijection,
and the associated one-to-one correspondence can be thought of as an identification.)
For simplicity of further notation, let us suppose that the $N$ ( $=$ number-of-vertices-choose-two) pairs
of vertices are arbitrarily ordered.

Define the set $\para:=\{ (p_1,p_2,\ldots,p_N, \varrho_1, \varrho_2, \ldots, \varrho_N):
p_1,p_2,\ldots,p_N, \varrho_1, \varrho_2, \ldots, \varrho_N \in [0,1]\}$.

\begin{definition} ({\bf Correlated Bernoulli Random Graph Model})
The parameter space for the correlated Bernoulli random graph model,
denoted $\Theta$, is any particular subset of $\para$, possibly a proper subset.
For each $(p_1,p_2,\ldots,p_N, \varrho_1, \varrho_2, \ldots, \varrho_N) \in \Theta$, the pair of random graphs
are described as follows. For each $i=1,2,\ldots,N$, the indicator
random variable $X_i$ for adjacency of the $i$th pair of vertices in the first graph and the indicator
random variable $Y_i$ for adjacency of the $i$th pair of vertices in the second graph are each
marginally distributed Bernoulli$(p_i)$, and the Pearson correlation coefficient of $X_i,Y_i$ is $\varrho_{i}$
(assuming that $p_i$ is not $0$ or $1$, in which case the value of $\varrho_i$ is irrelevant).
Other than these dependencies, the random variables $X_1,X_2,\ldots,X_N,Y_1,Y_2,\ldots,Y_N$ are independent.
\end{definition}

It is not hard to see that that the choice of these parameters uniquely specifies the joint
distribution of the two graphs (see Appendix A of \cite{TotCor}).
Indeed, we can sample from the distribution in the following manner. For each $i=1,2,\ldots,N$
independently, sample $X_i \sim \textup{Bernoulli}(p_i)$, then conditioned on the value $x_i$ of $X_i$,
sample $Y_i \sim \textup{Bernoulli}(\varrho_i x_i +(1-\varrho_i)p_i)$.
For all $i=1,2,\ldots, N$, define
\begin{eqnarray*}
q_{i,1} & := & p_i^2+\varrho_i p_i (1-p_i), \\
q_{i,0} & := & (1-p_i)^2+ \varrho_i p_i (1-p_i),\\
q_{i,\star} & := & (1-\varrho_i)p_i (1-p_i);
\end{eqnarray*}
these are, respectively, the probability that $X_i=Y_i=1$, the probability that $X_i=Y_i=0$, and the probability that
[$X_1=1$ and $Y_i=0$].

Let $\vx$ and $\vy$ denote the random vectors whose $i$th components are respectively $X_i$ and $Y_i$, for all
$i=1,2,\ldots,N$; thus, in effect, $\vx$ and $\vy$ are like the adjacency matrices representing the respective graphs.
Let $\x := \{ (x,y): x,y \in \{ 0,1 \}^N \}$ denote the sample space for the correlated
Bernoulli random graph model; in particular, $x$ and $y$ respectively are possible realizations of
the adjacency vectors $\vx$ and $\vy$.

 Note that if
$\varrho_1=\varrho_2=\cdots=\varrho_N=1$ then almost surely the two graphs are isomorphic, and if
$\varrho_1=\varrho_2=\cdots=\varrho_N=0$ then the two graphs are independent, meaning that the
collection of random variables
$X_1,X_2,\ldots,X_N,Y_1,Y_2,\ldots,Y_N$ is independent.

 Define
 $\paraoo:=\{ (p_1,p_2,\ldots,p_N,0,0, \ldots,0):p_1,p_2,\ldots,p_N \in \R \}$.
The parameter space $\Theta$ will be called {\it nondegenerate} if
$\Theta \cap \paraoo$ has an interior point, relative to $\paraoo$; i.e.~there exists
$z \in \Theta \cap \paraoo$ and real number $\epsilon>0$ such that $\Theta \cap \paraoo$ contains all points in $\paraoo$
that are less than $\epsilon$ distant~from~$z$.
Nondegeneracy of $\Theta$ will play a critical role here; it is an assumption explicitly required for
most of the theorems in this paper. Furthermore, when this condition is assumed, it is not
merely for ease of exposition or analysis; indeed, we will demonstrate that absence of this condition,
when it is assumed, can falsify the conclusions of the theorems that assume this condition.

\begin{remark}
\textup{ The results in this paper provide
machinery for
improved estimation in the context of the correlated Bernoulli random graph model,
which is a versatile and currently popular
random graph model utilized heavily in
the study of graph matching
and similar disciplines. Nonetheless,
the way the correlated Bernoulli random graph model is defined here---and the nature
of the results in this paper
and their proofs---render these results
also expressible more broadly
in terms of random correlated Bernoulli
vectors $X$,$Y$, for the pair of random
vectors $X$,$Y$ defined above, without
underlying graph structure.
Thus, in particular, with this broader perspective,
we do not need to restrict the
number of components $N$  (for $X$ and $Y$) to be ${n \choose 2}$, where $n$ is the number of vertices in an underlying graph. Indeed, we can consider $N$ as any positive integer.
} \end{remark}

\subsection{Important statistics and functions of the parameters \label{sec:stats}}

The most important statistic in this paper, the {\it disagreement vector} statistic $\HH$, is defined first.
This statistic is foundational for the first group of our results; in Theorem \ref{lem:complete} we
will show, under the nondegeneracy condition, that $\HH$ is complete and sufficient.

\begin{definition}
The (vector-valued) {\it disagreement vector} statistic
$\HH : \x \rightarrow \stuff$ is defined as follows: For all $(x,y) \in \x$, the vector $\HH(x,y) \in \stuff$
is such that, for each $i=1,2,\ldots, N$, the $i$th component of $\HH(x,y)$ is equal to $1$ if $x_i=y_i=1$,
is equal to $0$ if $x_i=y_i=0$, and is equal to $\star$ if $x_i \ne y_i$.  For all $h \in \stuff$,
the preimage of $h$, which is the set $\HH^{-1}(h)$, is denoted as $\x_h$, and is
called a {\it disagreement class}.
Note that $\x $ is partitioned into the disjoint union $\x =\bigcup_{h \in \stuff} \x_h $.
\end{definition}

The following definitions are key for the second group of our results.

The {\it Bernoulli parameter mean} $\mu$ and the {\it Bernoulli parameter variance} $\sigma^2$
are defined as
\begin{eqnarray*}
\mu:=\frac{1}{N}\sum_{i=1}^Np_i \  , \ \ \ \ \ \ \ \ \sigma^2:=\frac{1}{N}\sum_{i=1}^N(p_i-\mu)^2 \ .
\end{eqnarray*}

The {\it empirical density of $X$}, denoted $\dens_X$, the {\it empirical density of $Y$}, denoted $\dens_Y$,
and the {\it combined empirical density}, denoted $\dens_{X,Y}$, are statistics $\x \rightarrow \R$
that are respectively defined as
\begin{eqnarray*} \dens_X:=\frac{1}{N}\sum_{i=1}^NX_i \ , \ \ \ \ \
\dens_Y:=\frac{1}{N}\sum_{i=1}^NY_i \ , \ \ \  \ \
\dens_{X,Y}:=\frac{1}{2}(\dens_{X}+\dens_Y) \ \ .
\end{eqnarray*}
Clearly, all three of these statistics are unbiased estimators of the parameter $\mu$.
Then, we define the statistics $\densca, \denscu : \x \rightarrow \R$ as
\begin{eqnarray*}
\densca :=\frac{1}{N}\sum_{i=1}^NX_iY_i \ , \ \ \ \ \  \ \
\denscu := \dens_X+\dens_Y-\densca \ \ .
\end{eqnarray*}
Note that for all $(x,y) \in \x$, we have that  $\densca (x,y):=
\frac{| \{ i\ : x_i=1 \mbox{ and } y_i=1  \} |}{N}$, and we also have that
$\denscu (x,y):= \frac{| \{ i\ : x_i=1 \mbox{ or } y_i=1  \} |}{N} $.

Next, the {\it disagreement enumeration statistic} $\Delta: \x \rightarrow \R$ is
\begin{eqnarray*}
\Delta:=\sum_{i=1}^N(X_i-Y_i)^2  \ ;
\end{eqnarray*}
in particular, for all $(x,y) \in \x$, $\Delta(x,y)$
is the number of components at which $x$ and $y$ disagree. Clearly, we have, for all $\theta \in \Theta$, that
$\e (\Delta)=2\sum_{i=1}^N(1-\varrho_i)p_i(1-p_i)$. For all $h \in \stuff$ and  $(x,y) \in \x_h$,
we have $2^{\Delta(x,y)}=|\x_h|$.

The {\it heterogeneity correlation} $\varrho_H$ is a parameter defined by
\begin{eqnarray*}
\varrho_H:=\frac{\sigma^2}{\mu (1-\mu)} \ .
\end{eqnarray*}
In the case where $\mu$ is $0$ or $1$ then any convention may be adopted for defining $\varrho_H$ (but it must
be a value between $0$~and~$1$). It is not hard to show that {\bf a)}
it holds that $0 \leq \varrho_H \leq 1$, and
{\bf b)} it holds that $\varrho_H=1$ if and only if each of the $p_i$'s are either $0$ or $1$, and
{\bf c)} it holds that $\varrho_H=0$ if and only
if all of the $p_i$'s are equal to each other (of course, statements {\bf b)} and {\bf c)}
aren't to be applied to the case where $\mu$ is $0$ or $1$).

Define the {\it total correlation} parameter $\varrho_T$ as
\begin{eqnarray*}
\varrho_T:=1-\frac{\sum_{i=1}^N(1-\varrho_i)p_i(1-p_i)}{N \mu (1-\mu)} \ .
\end{eqnarray*}
In the case where $\mu$ is $0$ or $1$ then any convention may be adopted for defining $\varrho_T$ (but it must
be a value between $0$ and $1$). Note that in the case where all $\varrho_i$ are equal, say to the value $\varrho_E$,
then $(1-\varrho_T)=(1-\varrho_E)(1-\varrho_H)$. It is always the case that $0 \leq \varrho_T \leq 1$.
In \cite{TotCor} it was empirically demonstrated---for the correlated graphs in broad families within
our model---that graph matching complexity
as well as graph matchability are each functions of total correlation, hence the importance of total correlation.

The {\it alignment strength} $\st : \x \rightarrow \mathbb{R}$ is an important statistic defined as
\begin{eqnarray} \label{def:as}
\st :=1-\frac{ \Delta/N }{\dens_X \left ( 1-\dens_Y \right )+  \left ( 1-\dens_X \right ) \dens_Y} \ .
\end{eqnarray}
In the case that
$x$ and $y$ are both all zeros or both all ones then any convention may be adopted for defining $\st$ (but
it should be a value between $0$ and $1$).
The definition of alignment strength $\st$  arose in \cite{TotCor} in a natural way. Specifically, $1$ minus
the alignment strength is the ratio of disagreements between the two graphs ---through the natural alignment---
divided by the average number of disagreements over all vertex bijections between the two graphs;
see there for more details. In \cite{TotCor} it is proven under mild conditions that $\st$ is a strongly
consistent estimator of $\varrho_T$.

An equivalent formula for alignment strength is
\begin{eqnarray} \label{def:asother}
\st = \frac{\densca - \dens_X \dens_Y}{\dens_{X,Y} - \dens_X \dens_Y} \ ;
\end{eqnarray}
it follows immediately using the easily-derived identities mentioned later in
 Equations~(\ref{eqn:id1})~through~(\ref{eqn:id5}).

\section{The Results   \label{sec:results}}%%%%%%%%%%%%%

As mentioned in Section \ref{S:overview}, our main results, which will be listed in this section, can be divided
into two groups.

The first group of our results: In the context of correlated Bernoulli random
graphs,~we~begin with Theorem~\ref{lem:complete}, which asserts, under a nondegeneracy condition,
that the  disagreement vector statistic $\HH$ is complete and sufficient; using this,
given any estimator of a function of model parameters, we describe a way to refine (``balance")
the estimator, reducing the mean squared error. Indeed, under the nondegeneracy condition, given any unbiased estimator
of a function of model parameters, we
characterize all unbiased estimators (Theorem~\ref{thm:first})
and the UMVU estimator (Theorem~\ref{thm:second}).
The second group of our results: Theorems \ref{thm:hom}, \ref{thm:tot}, and \ref{thm:edg} show that
there are no unbiased estimators of various graph correlation measures, including
total correlation $\rT$; however, not only does balancing alignment strength improve alignment strength's mean
squared error in estimating $\rT$, but balancing numerator and denominator separately
is seen empirically to be a further improvement.

Our first result is the following theorem.

\begin{theorem} If the parameter space $\Theta$ is nondegenerate, then the disagreement vector
$\HH$ is a complete and sufficient statistic. \label{lem:complete}
\end{theorem}

Theorem \ref{lem:complete} is proved in Section \ref{sec:proofcomplete}.

Given any statistic $S: \x \rightarrow \R$, define the statistic
$\bS : \x \rightarrow \R$ as follows. For all $(x,y) \in \x$ and $h \in \stuff$ such that $(x,y) \in \x_h$,
define $\bS (x,y):= \frac{1}{|\x_h|}\sum_{(x',y')\in \x_h}S(x',y')$; in particular,
$\overline{S}$ is a constant function on $\x_h$. We say that $\bS$ is the {\it balanced}
variant of $S$; the {\it balancing} of $S$ means the substituting of $\bS$ in place of $S$
when performing an estimation or inference task.
If $S=\bS$ then we say $S$ is a {\it balanced} statistic. Of course, $S$ is balanced
if and only if $S$ is a constant function on $\x_h$ for each $h \in \stuff$. The following are our main results.

\begin{theorem} Suppose the parameter space $\Theta$ is nondegenerate, and a
statistic $S: \x \rightarrow \R$ is an unbiased estimator of $g(\theta)$, where $g: \Theta \rightarrow \R$.
Then a statistic $T: \x \rightarrow \R$ is an unbiased estimator of $g(\theta)$ if and only if
$\bT=\bS$. In particular, $\bS$ is an unbiased estimator of $g(\theta)$. \label{thm:first}
\end{theorem}

\begin{theorem} Suppose the parameter space $\Theta$ is nondegenerate, and a
statistic $S: \x \rightarrow \R$ is an unbiased estimator of $g(\theta)$, where $g: \Theta \rightarrow \R$.
Then there exists a UMVU estimator of $g(\theta)$ and, in fact, $\bS $ is the UMVU estimator of $g(\theta)$. \label{thm:second}
\end{theorem}

We prove Theorem \ref{thm:first} in Sections \ref{sec:firstreverse} and  \ref{sec:firstforward}, and
we prove Theorem \ref{thm:second} in Section \ref{sec:second}; it is essentially
a consequence of Theorem \ref{lem:complete}.
The key idea in proving
Theorem \ref{thm:second} is that $\bS $ is the composition of some function with $\HH$ ---a complete
and sufficient statistic by Theorem \ref{lem:complete}---and thus the Lehmann-Scheffe Theorem dictates that
$\bS$ is UMVU. Indeed, $\bS$ is the Rao-Blackwellization of $S$~conditioning~on~$\HH$.
Section \ref{sec:second} spells out the details.

(An excellent reference for Rao-Blackwell theory are the original papers \cite{Rao, Blackwell},
and an excellent reference for the Lehmann-Scheffe Theorem are the original papers \cite{LSone, LStwo}.)

\begin{remark} \textup{ Fundamentally, since balancing is Rao-Blackwellization, it is a
regularization technique that reduces an estimator's variance without changing its mean.
Indeed, suppose the parameter space $\Theta$ is nondegenerate, and a
statistic $S: \x \rightarrow \R$ is any estimator, whether biased or unbiased,
of $g(\theta)$. Then $\bS$ has the minimum possible variance
among all estimators with the same expected value as $S$.
(This, in turn, implies that $\bS$ has the minimum mean squared error
among all such estimators.) To see this, define
$\ell (\theta):=\e (S)$ for all $\theta \in \Theta$, and then
consider $S$, $\bS$, and $\ell$ in Theorems
\ref{thm:first} and \ref{thm:second}.
} \end{remark}

\begin{example} \label{ex:bal}
The disagreement enumeration statistic $\Delta: \x \rightarrow \R$ is clearly a balanced statistic, since
it is a constant function on each $\x_h$.
Hence, by Theorem \ref{thm:second}, when $\Theta$ is nondegenerate,
$\Delta$ is the UMVU estimator of its expected value
$\e (\Delta) = 2\sum_{i=1}^N(1-\varrho_i)p_i(1-p_i)$.
\end{example}

\begin{example} \label{ex:nobal}
When $N>1$, the statistic
$\dens_X \left ( 1-\dens_Y \right )+  \left ( 1-\dens_X \right ) \dens_Y$ is {\bf NOT} a balanced statistic;
indeed, consider
$h=[\star, \star, \ldots, \star]^T$, and consider $(x',y') \in \x_h$ such that
$x'$ is all zeros and $y'$ is all ones, and consider $(x'',y'') \in \x_h$ such that the
first $\lfloor \frac{N}{2} \rfloor$ entries of $x''$ are all zeros and of $y''$ are all ones, and
the last $\lceil \frac{N}{2} \rceil$ entries of $x''$ are all ones and of $y''$ are all zeros---the
statistic $\dens_X \left ( 1-\dens_Y \right )+  \left ( 1-\dens_X \right ) \dens_Y$
at $(x',y')$ has the value $1$, and at $(x'',y'')$ has a value approaching $\frac{1}{2}$,
hence the statistic is not constant on $\x_h$, hence is not balanced.
\end{example}

\begin{example} \label{ex:nonobal}
When $N>1$, the alignment strength statistic $\st$ is {\bf NOT} a balanced statistic.
This is because in Example \ref{ex:bal} we have that the numerator of $\st$ in Equation (\ref{def:as})
is balanced, and in Example \ref{ex:nobal} we have that the denominator of $\st$ in  Equation (\ref{def:as})
is not balanced; hence $\st$ is not a constant function on all $\x_h$, and is thus not balanced.
\end{example}

It is important to note that the claims in Theorems \ref{thm:first} and \ref{thm:second} may fail
without the assumption of nondegeneracy for the parameter space
$\Theta$, as highlighted in the next example.
(In particular, this points to the non-triviality of Theorems \ref{thm:first} and \ref{thm:second}.)

\begin{example} \label{ex:counter}
If the value of $\mu$ is known, hence fixed, then $\Theta$ is contained in a
particular hyperplane intersecting $\paraoo$, and $\Theta$ is degenerate; in this scenario,
we will show in Section~\ref{sec:need} that unbiasedness of estimators is
not characterized as described in Theorem \ref{thm:first} and, often,
there do not exist UMVU estimators for functions of the model parameters,
even when there exist unbiased estimators.
\end{example}

The following corollary is an immediate consequence of Theorem \ref{thm:second} and the fact
that sums and products of constant functions (on respective $\x_h$)
are constant (on respective $\x_h$).

\begin{corollary}  Suppose the parameter space $\Theta$ is nondegenerate, and a
statistic $S: \x \rightarrow \R$ is an unbiased estimator of $g(\theta)$,
where $g: \Theta \rightarrow \R$, and a statistic
$S': \x \rightarrow \R$ is an unbiased estimator of $g'(\theta)$, where $g': \Theta \rightarrow \R$.
Then, for any $a,b \in \R$, \ \  $a \bS+ b \bS'$ is UMVUE for $ag(\theta)+bg'(\theta)$. Indeed,
$a \bS + b \bS'$ is balanced, $\bS \cdot \bS'$ is balanced, and (if $\bS'$ is nonzero)
$\bS / \bS'$ is balanced.   \label{cor:vectorspace}
\end{corollary}

The next set of theorems are applications of the above theorems---and the methodologies of
their proofs---to unbiasedness and efficiency of statistics for estimating various
graph correlation parameters, particularly total correlation.

\begin{theorem} \label{thm:hom}
Suppose the parameter space $\Theta$ is nondegenerate, and $N>1$.
There does not exist an unbiased estimator of the heterogeneity correlation $\varrho_H$.
\end{theorem}

\begin{theorem} \label{thm:tot}
Suppose the parameter space $\Theta$ is nondegenerate, and $N>1$.
There does not exist an unbiased estimator of the total correlation $\varrho_T$.
\end{theorem}
We prove Theorems \ref{thm:hom} and \ref{thm:tot} in Section \ref{sec:tot}.

In the following negative result on estimating edge correlation,
besides the assumption of a nondegenerate parameter space, we have additional assumptions that all pairs of vertices
share the same edge correlation parameter (ie. $\Theta$ is restricted so that,
for all $\theta \in \Theta$, it holds that $\varrho_1=\varrho_2 \cdots = \varrho_N$),
and we also assume that this edge correlation~parameter~is~not~always~zero.~Specifically:

\begin{theorem} \label{thm:edg} Suppose that the following three conditions hold:\\
{\bf a)} The parameter space $\Theta$ is nondegenerate.\\
{\bf b)} The parameter space $\Theta$ is such that edge correlations are component-uniform, meaning that
there exists a function $\varrho_E : \Theta \rightarrow \R$ such that, for all
$\theta \in \Theta$ and all $i=1,2,\ldots,N$, \ $\varrho_i(\theta)=\varrho_E(\theta)$.\\
{\bf c)} The parameter space $\Theta$ is not a subset of $\paraoo$, i.e.~there exists $\theta \in \Theta$ such that $\varrho_E(\theta) \ne 0 $.\\
Then there does not exist an unbiased estimator of $\varrho_E$.
\end{theorem}
We prove Theorem \ref{thm:edg} in Section \ref{sec:edg}.

\begin{remark} \label{rem:quasi}
Suppose the parameter space $\Theta$ is nondegenerate. As mentioned in Example \ref{ex:nonobal},
alignment strength $\st = \frac{\densca - \dens_X \dens_Y}{\dens_{X,Y} - \dens_X \dens_Y}$ is {\bf NOT} balanced when $N>1$.
(So, the bias in estimating $\rT$ is the same for $\st$ as for $\overline{\st}$,
but the variance of $\overline{\st}$ is less than the variance of $\st$.) Next, define the  \textup{modified alignment strength}
$\st ' := \frac{ \ \ \ \overline{\densca - \dens_X \dens_Y} \ \ \ }{ \ \ \ \overline{\dens_{X,Y} - \dens_X \dens_Y} \ \ \ } $;
by Corollary \ref{cor:vectorspace}, $\st'$ is balanced.
We will empirically show in Section \ref{sec:empir} that $\st'$ is often superior to $\overline{\st}$ as an estimator of $\rT$.
Also, in Section \ref{sec:quasi} we will prove the following clean formulas:
\begin{eqnarray} \st' \ \ = \ \ \frac{\densca - \dens_X \dens_Y + \frac{1}{4} \left [ \frac{\Delta}{N^2} - (\dens_X-\dens_Y)^2 \right ]  }{
\dens_{X,Y} - \dens_X \dens_Y + \frac{1}{4} \left [ \frac{\Delta}{N^2} - (\dens_X-\dens_Y)^2 \right ]  } \ \ = \ \
\frac{\densca - \dens_{X,Y}^2 + \frac{\Delta}{4N^2}  }{ \dens_{X,Y} (1-\dens_{X,Y}) +\frac{\Delta}{4N^2}  }. \label{eqn:balnumdenom}
\end{eqnarray}
\end{remark}

When $x$ and $y$ are both all zeros or all ones then any convention for defining $\st '$ is acceptable,
provided that it is between $0$ and $1$. Note that $\st'$ can have less variance than $\overline{\st}$ (and in general
it does) without violating Theorem \ref{thm:second}, since the expected values of $\st'$ and $\overline{\st}$ can be  different.

Sometimes balancing a statistic---even at one sample space point---requires averaging an
exponential number of values. Remark \ref{rem:quasi} is notable for its simple expressions for the balanced
statistics comprising $\st'$, and in Section \ref{sec:empir} a linear time algorithm is given for
computing $\overline{\st}$.

The following result is proved in Section \ref{sec:cor}; it follows from Corollary \ref{cor:vectorspace}
and Remark \ref{rem:quasi}.

\begin{corollary} Suppose the parameter space $\Theta$ is nondegenerate and also suppose that $\Theta \subseteq \paraoo$.
Then the statistic $\dens_{X,Y} \left ( 1- \dens_{X,Y} \right ) -  \frac{1}{2N} \left ( 1- \frac{1}{2N} \right ) \Delta $ is
UMVUE for $\sigma^2$. \label{cor:cor}
\end{corollary}

These are the results in this paper, and they will be proven next in Section \ref{sec:pf}.

\section{Proof of the Results in Section \ref{sec:results}\label{sec:pf} }%%%%%%%%%%%%%

We begin by proving Theorem \ref{thm:first} and Theorem \ref{thm:second}; the proofs of Theorem \ref{lem:complete} and Theorem \ref{thm:edg}
will then be built on the methodology of the forward direction of the proof of Theorem \ref{thm:first}. The rest of the
results in Section~\ref{sec:results} will also be shown.

\subsection{Proof of the reverse direction of Theorem \ref{thm:first}   \label{sec:firstreverse}}

The reverse direction of Theorem \ref{thm:first} can be equivalently formulated in the following way.
Suppose the parameter space $\Theta$ is nondegenerate, and a
statistic $S: \x \rightarrow \R$ is an unbiased estimator of $g(\theta)$, where $g: \Theta \rightarrow \R$.
If the statistic $T: \x \rightarrow \R$ satisfies the condition that for all
$h \in \stuff \ \ \ \sum_{(x,y) \in \x_h} T(x,y) =  \sum_{(x,y) \in \x_h} S(x,y) $ then
$T$ is an unbiased estimator of $g(\theta)$.

Proving the reverse direction of Theorem \ref{thm:first} is quite straightforward.
For each $h \in \stuff$, the elements of $\x_h$ are equiprobable.
In particular, for all $\theta \in \Theta$ it holds that
\begin{eqnarray*}
\e (T) & = & \sum_{(x,y) \in \x} \Pry (x,y)  T(x,y) =  \sum_{h \in \stuff} \ \sum_{(x,y) \in \x_h} \Pry (x,y)  T(x,y) \\
& = & \sum_{h \in \stuff} \left [ \left ( \prod_{i:h_i=1}q_{i,1} \right )
\left ( \prod_{i:h_i=0}q_{i,0} \right ) \left ( \prod_{i:h_i=\star}q_{i,\star} \right )
\sum_{(x,y) \in \x_h} T(x,y) \right ] \\
& = & \sum_{h \in \stuff} \left [ \left ( \prod_{i:h_i=1}q_{i,1} \right )
\left ( \prod_{i:h_i=0}q_{i,0} \right ) \left ( \prod_{i:h_i=\star}q_{i,\star} \right )
\sum_{(x,y) \in \x_h} S(x,y) \right ] \\
&=&  \e (S) =  g(\theta) .
\end{eqnarray*}
Thus $T$ is an unbiased estimator of $g(\theta)$. $\qed$
\subsection{Proof of the forward direction of Theorem \ref{thm:first} and of Theorem \ref{lem:complete} \label{sec:firstforward}}

The proof of the forward direction of Theorem \ref{thm:first} involves notation that
is complex at first glance, and the core ideas may be challenging to follow when presented all at once in full generality. Our expositional
strategy is as follows. After proving the basic preliminary Lemma \ref{lem:polycoef} in Section~\ref{sec:prelim},
we proceed  to first prove  the forward direction of Theorem \ref{thm:first}
in the special cases where $N=1,2$ in Section~\ref{sec:special}, so that the notation, reasoning, and strategy are crystal clear,
and then in Section~\ref{sec:general} we prove the forward direction of Theorem \ref{thm:first} in
full and clear generality. Then we will show in Section \ref{sec:proofcomplete}
that the machinery of Section \ref{sec:general} proves Theorem \ref{lem:complete} (Completeness~of~$\HH$).

\subsubsection{Preliminaries \label{sec:prelim}}

We begin with a technical lemma, Lemma \ref{lem:polycoef}. Two polynomials in a single variable that
are equal as functions at infinitely many points are, by interpolation theory, equal algebraically
(meaning that the two polynomials have the same coefficients as each other). However,
for two polynomials in more than one variable, this may fail. For example, consider the
polynomial $p_1^2 -p_2$ and the zero polynomial, in the two variables $p_1$ and $p_2$. These
two polynomials agree as functions on a parabola, but they are not equal algebraically.
However, if two polynomials of any degree agree as functions on an open neighborhood
then they are equal algebraically. Formally:

\begin{lemma}
Suppose that $\Theta$ is nondegenerate, and $g, \tilde{g}: \Theta \cap \paraoo \rightarrow \R$ are two
polynomials in the variables $p_1,p_2,\ldots, p_N$ such that for
all $\theta \in \Theta \cap \paraoo$ it holds that $g(\theta)=\tilde{g}(\theta)$. Then the  coefficients
of the polynomial $g$ are identical to the respective coefficients of the polynomial $\tilde{g}$. \label{lem:polycoef}
\end{lemma}

The proof of Lemma \ref{lem:polycoef} is a straightforward induction on the maximum degree of
the polynomials $g$ and $\tilde{g}$, considering sequential partial derivatives. An equivalent
formulation of Lemma \ref{lem:polycoef} can be found in the classical textbook {\it Algebra}
of Serge Lang \cite{Lang}, Chapter IV, Corollary 1.6.

%the base case statement---about constant polynomials---is trivially true.
%Assuming the lemma statement is true for any two polynomials (equal as functions on $\Theta \cap \paraoo$)
%with degree at most $k$, consider any $g$ and $\tilde{g}$ (equal as functions on $\Theta \cap \paraoo$) with degree at most $k+1$.
%By the nondegeneracy of $\Theta$, we have that
%$\Theta \cap \paraoo$ has an interior point, relative to $\paraoo$;
%at that point, $k+1$ sequential applications of
%partial derivatives of $g$ and $\tilde{g}$ (with respect to various combinations of $p_1,p_2,\ldots,p_N$)
%yields the polynomial coefficients of the monomials of degree $k+1$ for $g$ and $\tilde{g}$,
%and such coefficients for $g$ are the same as for $\tilde{g}$. Subtracting these monomials from $g$ and $\tilde{g}$ we obtain two
%polynomials (equal as functions on $\Theta \cap \paraoo$)
%with degree at most $k$, and by the induction hypothesis the remaining monomials of $g$ are identical to those of $\tilde{g}$.
%The result follows by induction. $\qed$\\

\subsubsection{Proof of the forward direction of Theorem \ref{thm:first} for the particular cases where $N=1,2$, by way of
illustration for the general case \label{sec:special}}

The forward direction of Theorem \ref{thm:first} can be formulated as follows.
Suppose the parameter space $\Theta$ is nondegenerate, and the two
statistics $S,T: \x \rightarrow \R$ are each unbiased estimators of $g(\theta)$, where $g: \Theta \rightarrow \R$.
Then for all $h \in \stuff$ it holds that $\sum_{(x,y) \in \x_h} T(x,y) =  \sum_{(x,y) \in \x_h} S(x,y) $.

To best illustrate, we begin with a proof for the case where $N=1$. Taking the expectation for parameters $\theta \in \Theta \cap \paraoo$,
we see that $\e (S)= (1-p_1)^2 S(0,0) + p_1(1-p_1) \big (S(0,1)+S(1,0) \big ) + p_1^2S(1,1)$. In particular,
$g$ needs to be a quadratic polynomial in the single variable $p_1$ on  $\Theta \cap \paraoo$, say $g(p_1):=g^{(0)}p_1^0+g^{(1)}p_1^1+g^{(2)}p_1^2$
where $g^{(0)}$, $g^{(1)}$, and $g^{(2)}$ are fixed coefficients.
By the nondegeneracy of $\Theta$ and Lemma \ref{lem:polycoef}, we can uniquely
represent polynomials as vectors with respective entries
being the coefficients of $p_1^0$, $p_1^1$, and $p_1^2$, respectively.
Thus $(1-p_1)^2$ is represented as  $\bigl [ \begin{smallmatrix} 1 \\ -2 \\ 1  \end{smallmatrix} \bigr ]$,
and $p_1(1-p_1)$ is represented as $\bigl [ \begin{smallmatrix} 0 \\ 1 \\ -1  \end{smallmatrix} \bigr ]$,
and $p_1^2$ is represented as $\bigl [ \begin{smallmatrix} 0 \\ 0 \\ 1  \end{smallmatrix} \bigr ]$,
and $g$ is represented as $\bigl [ \begin{smallmatrix} g^{(0)} \\ g^{(1)} \\ g^{(2)}  \end{smallmatrix} \bigr ]$.
In particular, $S$ being an unbiased estimator for $g$ on $\Theta \cap \paraoo$ means that
\[ \left [  \begin{array}{rrr} 1 & 0 & 0 \\ -2 & 1 & 0 \\ 1 & -1 & \ \ 1    \end{array}  \right ]
   \left [   \begin{array}{c} S(0,0) \\ S(0,1)+S(1,0) \\ S(1,1)     \end{array}  \right ] =
    \left [   \begin{array}{c} g^{(0)} \\ g^{(1)} \\ g^{(2)}   \end{array}  \right ]
\]
Denote the left hand side matrix as $A$, that is
$A= \bigl [ \begin{smallmatrix} 1 & 0 & 0 \\ -2 & 1 & 0 \\ 1 & -1 & 1  \end{smallmatrix} \bigr ]$;
since $A$ is invertible, and $T$ has to satisfy the above equation as well, we therefore have that
 $\bigl [ \begin{smallmatrix} S(0,0) \\ S(0,1)+S(1,0) \\ S(1,1) \end{smallmatrix} \bigr ]=
 \bigl [ \begin{smallmatrix} T(0,0) \\ T(0,1)+T(1,0) \\ T(1,1) \end{smallmatrix} \bigr ]$, which precisely
 says that for all $h \in \stuff$ it holds that $\sum_{(x,y) \in \x_h} T(x,y) =  \sum_{(x,y) \in \x_h} S(x,y) $,
 and the case where $N=1$ is proven.

By further way of illustration, we next prove the case where $N=2$.  Taking the expectation for parameters $\theta \in \Theta \cap \paraoo$, we have
\begin{eqnarray*}
\e (S) & = & (1-p_1)^2(1-p_2)^2 \left [  S( \bigl [ \begin{smallmatrix} 0 \\ 0   \end{smallmatrix} \bigr ],
\bigl [ \begin{smallmatrix} 0 \\ 0   \end{smallmatrix} \bigr ] )   \right ]
+ (1-p_1)^2p_2(1-p_2) \Big [  S( \bigl [ \begin{smallmatrix} 0 \\ 0   \end{smallmatrix} \bigr ],
\bigl [ \begin{smallmatrix} 0 \\ 1   \end{smallmatrix} \bigr ] )+S( \bigl [ \begin{smallmatrix} 0 \\ 1   \end{smallmatrix} \bigr ],
\bigl [ \begin{smallmatrix} 0 \\ 0   \end{smallmatrix} \bigr ] )   \Big ] \ \ \ + \ \ \ \cdots  \\
& & \ \ \  \cdots \ \ \ + \ \ \ p_1(1-p_1)p_2(1-p_2) \Big [
 S( \bigl [ \begin{smallmatrix} 0 \\ 0   \end{smallmatrix} \bigr ],
\bigl [ \begin{smallmatrix} 1 \\ 1   \end{smallmatrix} \bigr ] )+
 S( \bigl [ \begin{smallmatrix} 0 \\ 1   \end{smallmatrix} \bigr ],
\bigl [ \begin{smallmatrix} 1 \\ 0   \end{smallmatrix} \bigr ] )+
 S( \bigl [ \begin{smallmatrix} 1 \\ 0   \end{smallmatrix} \bigr ],
\bigl [ \begin{smallmatrix} 0 \\ 1   \end{smallmatrix} \bigr ] )+
 S( \bigl [ \begin{smallmatrix} 1 \\ 1   \end{smallmatrix} \bigr ],
\bigl [ \begin{smallmatrix} 0 \\ 0   \end{smallmatrix} \bigr ] )
      \Big ]           \ \ \ + \ \ \ \cdots       \\
& = &  (1-p_1)^2(1-p_2)^2 \sum_{(x,y) \in \x_{\bigl [ \begin{smallmatrix} 0 \\ 0   \end{smallmatrix} \bigr ]}} S(x,y) \ \ \ + \ \ \
(1-p_1)^2p_2(1-p_2) \sum_{(x,y) \in \x_{\bigl [ \begin{smallmatrix} 0 \\ \star   \end{smallmatrix} \bigr ]}} S(x,y) \\
& &  + \ \ \  (1-p_1)^2p_2^2 \sum_{(x,y) \in \x_{\bigl [ \begin{smallmatrix} 0 \\ 1   \end{smallmatrix} \bigr ]}} S(x,y) \ \ \ + \ \ \
p_1(1-p_1)(1-p_2)^2 \sum_{(x,y) \in \x_{\bigl [ \begin{smallmatrix} \star \\ 0   \end{smallmatrix} \bigr ]}} S(x,y) \\
& &  + \ \ \  p_1(1-p_1)p_2(1-p_2) \sum_{(x,y) \in \x_{\bigl [ \begin{smallmatrix} \star \\ \star   \end{smallmatrix} \bigr ]}} S(x,y) \ \ \ + \ \ \
p_1(1-p_1)p_2^2 \sum_{(x,y) \in \x_{\bigl [ \begin{smallmatrix} \star \\ 1   \end{smallmatrix} \bigr ]}} S(x,y) \\
& &  + \ \ \  p_1^2(1-p_2)^2 \sum_{(x,y) \in \x_{\bigl [ \begin{smallmatrix} 1 \\ 0   \end{smallmatrix} \bigr ]}} S(x,y) \ \ \ + \ \ \
p_1^2p_2(1-p_2) \sum_{(x,y) \in \x_{\bigl [ \begin{smallmatrix} 1 \\ \star   \end{smallmatrix} \bigr ]}} S(x,y) \\
& &  + \ \ \  p_1^2p_2^2 \sum_{(x,y) \in \x_{\bigl [ \begin{smallmatrix} 1 \\ 1   \end{smallmatrix} \bigr ]}} S(x,y) .
\end{eqnarray*}
Note in particular that $g$ would have to be a polynomial in the two variables $p_1,p_2$,
with its monomials each consisting of a constant, denoted $g^{(k_1,k_2)}$,  times $p_1^{k_1}p_2^{k_2}$, where $k_1,k_2 \in \{0,1,2\}$.
By the nondegeneracy of $\Theta$ and Lemma \ref{lem:polycoef}, $g$ can be uniquely represented by the vector of
coefficients ordered lexicographically (i.e.~dictionary order) by superscript:
$[g^{(0,0)},g^{(0,1)}, g^{(0,2)},g^{(1,0)}, g^{(1,1)},g^{(1,2)}, g^{(2,0)},g^{(2,1)}, g^{(2,2)}]^T$.
Indeed, all other polynomials with monomials each consisting of a constant times $p_1^{k_1}p_2^{k_2}$, where $k_1,k_2 \in \{0,1,2\}$, will also
be similarly represented (uniquely) by the vector of coefficients ordered lexicographically by superscript.
For example, in the matrix on the left hand side below, the columns respectively are the vectors
(of lexicographically ordered monomial coefficients)  representing the respective polynomials
$(1-p_1)^2(1-p_2)^2$, \ $(1-p_1)^2p_2(1-p_2)$, \ $(1-p_1)^2p_2^2$, \dots, which are the respective probabilities
of $(x,y) \in \x_h$ for $h$'s lexicographically ordered (``$\star$"~has~the~value~$\frac{1}{2}$) as:
$\bigl [ \begin{smallmatrix} 0 \\ 0   \end{smallmatrix} \bigr ]$,
$\bigl [ \begin{smallmatrix} 0 \\ \star   \end{smallmatrix} \bigr ]$,
$\bigl [ \begin{smallmatrix} 0 \\ 1   \end{smallmatrix} \bigr ]$,
$\bigl [ \begin{smallmatrix} \star \\ 0   \end{smallmatrix} \bigr ]$,
$\bigl [ \begin{smallmatrix} \star \\ \star   \end{smallmatrix} \bigr ]$,
$\bigl [ \begin{smallmatrix} \star \\ 1   \end{smallmatrix} \bigr ]$,
$\bigl [ \begin{smallmatrix} 1 \\ 0   \end{smallmatrix} \bigr ]$,
$\bigl [ \begin{smallmatrix} 1 \\ \star   \end{smallmatrix} \bigr ]$,
$\bigl [ \begin{smallmatrix} 1 \\ 1   \end{smallmatrix} \bigr ]$.
Now, $S$ being an unbiased estimator of $g$ on $\Theta \cap \paraoo$ means precisely that the following linear system holds:
\begin{eqnarray} \label{eqn:specialcase}
\left [  \begin{array}{rrr|rrr|rrr}
1  & 0 & 0   &    0 & 0 & 0   &   0 & 0 & \ \ 0 \\
-2 & 1 & 0   &    0 & 0 & 0   &  0 & 0 & 0 \\
1  & -1 & 1  &    0 & 0 & 0   &  0 & 0 & 0 \\ \hline
-2 &  0 & 0  &    1 & 0 & 0   &  0 & 0 & 0 \\
4  & -2 & 0  &   -2 & 1 & 0   &  0 & 0 & 0 \\
-2 &  2 & -2 &    1 & -1 & 1  &  0 & 0 & 0 \\ \hline
1  & 0 & 0   &   -1 & 0 & 0   &  1 & 0 & 0 \\
-2 & 1 & 0   &    2 & -1 & 0   & -2 & 1 & 0 \\
1  & -1 & 1  &   -1 & 1 & -1   &  1 & -1 & 1 \\
\end{array}  \right ]
\left [   \begin{array}{c}
\sum_{(x,y) \in \x_{[0,0    ]^T}} S(x,y) \\
\sum_{(x,y) \in \x_{[0,\star]^T}} S(x,y) \\
\sum_{(x,y) \in \x_{[0, 1   ]^T}} S(x,y) \\ \hline
\sum_{(x,y) \in \x_{[\star, 0]^T}} S(x,y) \\
\sum_{(x,y) \in \x_{[\star,\star]^T}} S(x,y) \\
\sum_{(x,y) \in \x_{[\star, 1]^T}} S(x,y) \\ \hline
\sum_{(x,y) \in \x_{[1,0 ]^T}} S(x,y) \\
\sum_{(x,y) \in \x_{[1,\star]^T}} S(x,y) \\
\sum_{(x,y) \in \x_{[1,1]^T}} S(x,y)
\end{array} \right ] =
\left [   \begin{array}{c}
g^{(0,0)} \\ g^{(0,1)} \\ g^{(0,2)} \\ \hline
g^{(1,0)} \\ g^{(1,1)} \\ g^{(1,2)} \\ \hline
g^{(2,0)} \\ g^{(2,1)} \\ g^{(2,2)}
\end{array}  \right ]
\end{eqnarray}
Observe that the left-hand-side matrix above is the Kronecker product
$A \otimes A$, where $A$ is the lower triangular matrix with diagonals all ones mentioned in the proof of the case where $N=1$.
Note that $A \otimes A$ is thus lower triangular with diagonals all ones, thus has nonzero determinant and is invertible. Since $T$ solves the same
 linear system (above, Equation \ref{eqn:specialcase})
 as $S$ does, we conclude ---by multiplying both sides of the equation above by the
inverse of $A \otimes A$--- that for all $h \in \stuff$ it holds that
$\sum_{(x,y) \in \x_h} T(x,y) =  \sum_{(x,y) \in \x_h} S(x,y) $,
 and the case where $N=2$ is now also proven.

\subsubsection{Proof of the forward direction of Theorem \ref{thm:first}, the general case \label{sec:general}}

With the proofs of the cases $N=1,2$ as illustration, we now prove the result for arbitrary~$N$.
Let $\otimes^NA$ denote the $N$-fold Kronecker product $A \otimes A \otimes \cdots \otimes A$.
Next, let $\overrightarrow{ S/{\mathcal H}}$ denote the vector whose components are respectively
$\sum_{(x,y) \in \x_h} S(x,y)$ for each  of the $h \in \stuff$, ordered lexicographically according to $h$.
Restricting to parameters  $\theta \in \Theta \cap \paraoo$, $g$ is a polynomial in the variables $p_1,p_2,\ldots,p_N$,
with its monomials each consisting of a constant, denoted $g^{(k_1,k_2,\ldots,k_N)}$,
times $p_1^{k_1}p_2^{k_2}\cdots p_N^{k_N}$, where $k_1,k_2,\ldots,k_N \in \{0,1,2\}$.
By the nondegeneracy~of~$\Theta$ and Lemma \ref{lem:polycoef}, we have that $g$ can be uniquely represented by the column vector of monomial
coefficients ordered lexicographically by the powers of $p_1,p_2,\ldots,p_N$;  denote this vector $\vec{g}$.

We claim that $S$ being an unbiased estimator for $g(\theta)$ on $\Theta \cap \paraoo$ means precisely that $S$ satisfies
the linear system
\begin{eqnarray} \label{eqn:first}
[\otimes^N A] \  \cdot  \ \overrightarrow{ S/{\mathcal H}} \ = \ \vec{g}.
\end{eqnarray}
This can be verified directly by noting that for each $h \in \stuff$ and $(x,y) \in \x_h$, the probability of $(x,y)$ is
given by (and is simplified with elementary algebra)
\begin{eqnarray} \label{eqn:second}
\prod_{i=1}^N \left \{ \begin{array}{cl}
(1-p_i)^2 & \mbox{if } h_i=0 \\
p_i(1-p_i) & \mbox{if } h_i =\star \\
p_i^2 & \mbox{if } h_i =1
\end{array} \right \} \
=
\sum_{(k_1,k_2,\ldots,k_N)\in \{0,1,2\}^N}
\left ( \prod_{j=1}^N A_{k_j+1 \, , \, 2\cdot h_j+1} \right ) \cdot
\left ( \prod_{j=1}^N p_j^{k_j} \right ) \ ,
\end{eqnarray}
where, in the subscript of $A_{k_j+1 \, , \, 2\cdot h_j+1}$, ``$\star$" has the value $\frac{1}{2}$, meaning that
when $h_j$ is $\star$ then $2 \cdot ``\star" +1$ is defined to be $2$.
With nondegeneracy of $\Theta$ and Lemma~\ref{lem:polycoef}, we have uniqueness of polynomial coefficients,
and Equation (\ref{eqn:second}) directly yields Equation~(\ref{eqn:first}).

Since $\otimes ^N A$ is a lower triangular matrix with all diagonals being ones, it is an invertible matrix.
Now, $T$ has to satisfy Equation (\ref{eqn:first}) as well; multiplying both sides of the equation by the inverse of
$\otimes ^N A$ yields that $\overrightarrow{ S/{\mathcal H}} = \overrightarrow{ T/{\mathcal H}}$, which
precisely says that for all $h \in \stuff$ it holds that
$\sum_{(x,y) \in \x_h} T(x,y) =  \sum_{(x,y) \in \x_h} S(x,y) $,
and the forward direction of Theorem \ref{thm:first} is now proved. $\qed$

\subsubsection{Proof of Theorem \ref{lem:complete}  \label{sec:proofcomplete}}

Theorem \ref{lem:complete} states that
if the parameter space $\Theta$ is nondegenerate, then
$\HH$ is a complete and sufficient statistic. Sufficiency of $\HH$ is immediate,
since, conditioning on $\HH$  being any $h \in \stuff$,
we have that $(X,Y)$ is distributed discrete uniform on $\x_h$ (its support).
Also note that the probability of any sample point is a function of $\HH$ and the
model parameters, which implies that $\HH$ is sufficient.

We now use the machinery of the previous Section \ref{sec:general} to
prove the rest of Theorem \ref{lem:complete}; that if the parameter space $\Theta$ is nondegenerate,
then $\HH$ is a complete statistic. Completeness of $\HH$ here means that, for
any $f: \stuff \rightarrow \R$, if $\e (f(\HH))=0$ for all $\theta \in \Theta$ then $f$ is~the~zero~function.

For any $f: \stuff \rightarrow \R$, let $\vecf$ denote a column vector whose respective
entries are $|\x_h| \cdot f(h)$ for the respective $h \in \stuff$ ordered lexicographically. (I.e., the
first component of $\vecf$ is $|\x_h| \cdot f(h)$ for the $h$ that is all zeros, the second component
of $\vecf$ is $|\x_h| \cdot f(h)$ for the $h$ that is all zeros except last entry is $\star$, etc.)  By the nondegeneracy
of $\Theta$ and Equation (\ref{eqn:first}), we have that $\e (f(\HH))=0$ for all $\theta \in \Theta \cap \paraoo$ would mean that
$[\otimes^N A] \  \cdot  \ \vecf = \ \vec{0}$; by the invertibility of $[\otimes^N A]$, we would
have that $\vecf = \vec{0}$, hence $f$ is the zero function. This proves the completeness of $\HH$. $\qed$

Note that proof of the forward direction of Theorem \ref{thm:first} is
based on the injectivity of $[\otimes^N A]$ as a function, and proof of the completeness of $\HH$ in
Theorem \ref{lem:complete} is based on  $[\otimes^N A]$ having a trivial nullspace, so
these two results are equivalent.

\subsection{Proof of Theorem \ref{thm:second} \label{sec:second}}

Theorem \ref{thm:second} states that if the parameter space $\Theta$ is nondegenerate, and a
statistic $S: \x \rightarrow \R$ is an unbiased estimator of $g(\theta)$, for $g: \Theta \rightarrow \R$,
then there exists a UMVU estimator of $g(\theta)$ and, in fact, the balanced statistic $\bS $ is the UMVU estimator of $g(\theta)$.

We prove this as follows, under the assumption that $\Theta$ is nondegenerate.
Recall first that $\bS$ is an unbiased estimator of $g(\theta)$ by Theorem \ref{thm:first}.
Now, note that $\bS$ is a constant function on $\x_h$ for each $h \in \stuff$; thus there exists a
function $\Phi : \stuff \rightarrow \R$ such that $\bS$ is the function composition $\Phi \circ \HH$.
Since $\HH$ is a complete and sufficient statistic by Theorem \ref{lem:complete}, the Lehmann-Scheffe Theorem
asserts that the composition  $\bS=\Phi \circ \HH$ is UMVUE for $g(\theta)$. $\qed$

It is structurally interesting to note that $\bS$ is the Rao-Blackwellization of $S$ when conditioning
on the complete and sufficient statistic $\HH$. Indeed, for any $h \in \stuff$, the Rao Blackwellization
$\e ( S | \HH =h ) $ is precisely
the mean of the values of $S$ on $\x_h$, which is precisely the statistic $\bS$.

\subsubsection{Another proof of Theorem \ref{thm:second}, by first-principles \label{sec:elem}}

In this section, we mention a ``first-principles" proof of Theorem \ref{thm:second}, besides
the Lehmann-Scheffe proof methodology in the previous Section \ref{sec:second}.

For any statistic $T: \x \rightarrow \R$ which is an unbiased estimator of $g(\theta)$,
we compute the variance of $T$, for any particular $\theta \in \Theta$, as:
\begin{eqnarray*}
\textup{Var} (T) & = & \sum_{(x,y) \in \x} \Pry (x,y) \big ( T(x,y) - g(\theta) \big )^2 \\
& = & \sum_{h \in \stuff} \ \  \sum_{(x,y) \in \x_h} \Pry (x,y) \big ( T(x,y) - g(\theta) \big )^2\\
& = & \sum_{h \in \stuff} \left [ \left ( \prod_{i:h_i=1}q_{i,1} \right )
\left ( \prod_{i:h_i=0}q_{i,0} \right ) \left ( \prod_{i:h_i=\star}q_{i,\star} \right )
\sum_{(x,y) \in \x_h} \big ( T(x,y) - g(\theta) \big )^2.
\right ] \end{eqnarray*}
In particular, Var$(T)$ can be minimized over such unbiased $T$ by, for all $h \in \stuff$,
minimizing $\sum_{(x,y) \in \x_h} \big ( T(x,y) - g(\theta) \big )^2$ subject to the constraint that
$\sum_{(x,y) \in \x_h} T(x,y) =  \sum_{(x,y) \in \x_h} S(x,y)$; this is because of Theorem \ref{thm:first}.
Treating the $T(x,y)$ as variables, this convex optimization
problem has a global minimizer, when the objective gradient is equivalued (by the KKT conditions) hence the minimum variance
is achieved when $T=\bS$, independent of $\theta \in \Theta$, and Theorem \ref{thm:second} is shown.~$\qed$

\subsection{Proof of Theorem \ref{thm:edg}  \label{sec:edg}}

\noindent Theorem \ref{thm:edg} states that if the following three conditions hold:\\
{\bf a)} The parameter space $\Theta$ is nondegenerate.\\
{\bf b)} The parameter space $\Theta$ is such that the edge correlations are component-uniform, meaning that
there exist a function $\varrho_E : \Theta \rightarrow \R$ such that, for all
$\theta \in \Theta$ and all $i=1,2,\ldots,N$, \ $\varrho_i(\theta)=\varrho_E(\theta)$.\\
{\bf c)} The parameter space $\Theta$ is not a subset of $\paraoo$, i.e.~there exists $\theta \in \Theta$ such that $\varrho_E(\theta) \ne 0 $.\\
 Then there does not exist an unbiased estimator of $\varrho_E(\theta)$.\\

Suppose, by way of contradiction, that statistic $S: \x \rightarrow \R$ is an unbiased estimator of $\varrho_E(\theta)$.
For $\theta \in
\Theta \cap \paraoo$, where $\varrho_E \equiv 0$,
we have by Equation (\ref{eqn:first}) that $[\otimes^N A] \  \cdot  \ \overrightarrow{ S/{\mathcal H}} \ = \ \vec{0}$, since $\Theta$ is nondegenerate.
By the invertibility of $\otimes^N A$ we thus have that $\overrightarrow{ S/{\mathcal H}} \ = \ \vec{0}$,
which, by the reasoning in Section \ref{sec:firstreverse}, implies that $\e (S) =0$ for all $\theta \in \Theta$, which is a contradiction because there exists $\theta \in \Theta$
where $\e (S) = \varrho_E(\theta) \ne 0$. $\qed$

\subsection{Proof of Theorems \ref{thm:hom} and \ref{thm:tot} \label{sec:tot} }
Theorems \ref{thm:hom} and \ref{thm:tot} state that if
the parameter space $\Theta$ is nondegenerate and $N>1$ then there does not exist an unbiased estimator of the
heterogeneity correlation $\varrho_H$ nor of the total correlation $\varrho_T$.

We will just focus on $\Theta \cap \paraoo$; on this set it is easy to see that $\varrho_T=\varrho_H$. Thus,
by the development in Section \ref{sec:firstforward} and the nondegeneracy of $\Theta$, we will have
proved Theorems \ref{thm:hom} and \ref{thm:tot} if we show that, on $\Theta \cap \paraoo$, \ \
$\varrho_H:=\frac{\sigma^2}{\mu (1-\mu)}$ is not a polynomial in the variables $p_1,p_2,\ldots,p_N$.
By way of contradiction, suppose that,  on $\Theta \cap \paraoo$, \ \
$\varrho_H$ is a polynomial in the variables $p_1,p_2,\ldots,p_N$.

Let $(\tilde{p_1},\tilde{p_2},\ldots,\tilde{p_N},0,0,\dots,0)$ be an interior point of $\Theta \cap \paraoo$, relative to $\paraoo$;
such a point exists by the nondegeneracy of $\Theta$.
%%sew
Consider fixing the values of $p_i$ to be $\tilde{p_i}$ for each $i=2,3,\ldots,N$, and varying only $p_1$. This
results in $\varrho_H$, $\sigma^2$, and $\mu$ being polynomials in a single variable. Denote this variable
by $t$ instead of $p_1$ for ease of notation, and these respective polynomials are thus denoted  $\varrho_H(t)$, $\sigma^2(t)$, and $\mu (t)$.
%Let $\varrho_H(t)$, $\sigma^2(t)$, and $\mu (t)$ respectively
%denote the polynomials in the single variable $t$, where $t$ is substituted in place of $p_1$ and, for each of $i=2,3,\ldots,N$, the value of
%$p_i$ is fixed to be $\tilde{p_i}$, for $\varrho_H$, $\sigma^2$, and $\mu$ respectively.
%%sew
Let ${\mathcal I}$ be a real, open interval
containing $\tilde{p_1}$, such that for all $t \in {\mathcal I}$ we have
$(t,\tilde{p_2},\ldots,\tilde{p_N},0,0,\dots,0) \in \Theta \cap \paraoo$; a nontrivial such ${\mathcal I}$ exists by the
nondegeneracy of $\Theta$.

Using basic algebra, $\sigma^2(t)$ is quadratic in $t$, and the coefficient of $t^2$ in $\sigma^2(t)$
is $\left ( \frac{1}{N} -\frac{1}{N^2} \right )$, and $\mu(t)(1-\mu(t))$ is quadratic in $t$, and the coefficient of $t^2$
in $\mu(t)(1-\mu(t))$ is $-\frac{1}{N^2}$. Now, by definition, $\sigma^2(t)= \mu(t)(1-\mu(t))\varrho_H(t)$, and the coefficients
of the respective powers of $t$ on the left hand side are respectively equal to the coefficients of the powers of $t$ on the right
hand side, since ${\mathcal I}$ is an interval (and invoking polynomial interpolation theory). This implies that polynomial
$\varrho_H(t)$ can't have positive degree, and thus is constant. However, this constant is nonnegative (indeed,
it has been pointed out in Section \ref{sec:stats} that $0 \leq \varrho_H \leq 1$), which means that the coefficient of
$t^2$ in (the left hand side) $\sigma^2(t)$ is positive, but the coefficient of $t^2$ in (the right hand side)
 $\mu(t)(1-\mu(t))\varrho_H(t)$ is nonnegative times negative, which is nonpositive. By the contradiction, we have thus proved Theorems~\ref{thm:hom}~and~\ref{thm:tot}.~$\qed$

%\subsubsection{Cramer-Rao Lower Bound may be meaningless to establish UMVU}
% ************I COMMENTED THIS SECTION OUT, BUT IT WAS A FINE PART OF PREVIOUS VERSION OF MANUSCRIPT
%It is worthy of note that our proof of Theorem \ref{thm:second} characterizing UMVU
%estimators (when $\Theta$ is nondgenerate) does not make
%use of the Cramer-Rao lower bound. Indeed, the Cramer-Rao lower bound can fail to
%give a tight bound on the variance; for example, consider the disagreement statistic $\Delta$.
%For simplicity, suppose that $\Theta = \parao$ and, again for simplicity just over here, all $\varrho_i$'s
%are removed from the model, and $\e (\Delta)= 2\sum_{i=1}^Np_i(1-p_i)$.
%The Fisher information matrix is straightforward; its inverse is the diagonal matrix
%$\textup{diag}(\frac{p_1(1-p_1)}{2},\frac{p_2(1-p_2)}{2},\ldots, \frac{p_N(1-p_N)}{2})$.
%However, the gradient of $\e (\Delta)$ is  $[2-4p_1, 2-4p_2, \ldots, 2-4p_N]^T$, which is zero
%when $p_1=p_2=\cdots=p_N=\frac{1}{2}$, at which the Cramer-Rao lower bound is meaningless.
%Also in this case, the gradient of total correlation (note: $\varrho_T=\varrho_H=\frac{\sigma^2}{\mu (1-\mu)}$ when $\Theta = \parao$)
%is zero as well when $p_1=p_2=\cdots=p_N=\frac{1}{2}$.
%(This is seen by taking partial derivatives of $\sigma^2$ and also $\mu (1-\mu)$ with
%respect to the $p_i$; these partial derivatives are all zero. Then apply the quotient rule for derivatives.) Again,
%the Cramer-Rao lower bound is meaningless in this case.

\subsection{Proof of Equation (\ref{eqn:balnumdenom}) in Remark \ref{rem:quasi}  \label{sec:quasi}}

Recall that  $\st = \frac{\densca - \dens_X \dens_Y}{\dens_{X,Y} - \dens_X \dens_Y}$, and also recall the definition of the
modified alignment strength
$\st ' := \frac{ \ \ \ \overline{\densca - \dens_X \dens_Y} \ \ \  }{  \ \ \ \overline{\dens_{X,Y} - \dens_X \dens_Y} \ \ \  }$.
The main goal of this section is to prove Equation (\ref{eqn:balnumdenom}) in Remark \ref{rem:quasi}; namely, we show that
\begin{eqnarray}
\st' \ \ = \ \ \frac{\densca - \dens_X \dens_Y + \frac{1}{4} \left [ \frac{\Delta}{N^2} - (\dens_X-\dens_Y)^2 \right ]  }{
\dens_{X,Y} - \dens_X \dens_Y + \frac{1}{4} \left [ \frac{\Delta}{N^2} - (\dens_X-\dens_Y)^2 \right ]  } \ \ = \ \
\frac{\densca - \dens_{X,Y}^2 + \frac{\Delta}{4N^2}  }{ \dens_{X,Y} (1-\dens_{X,Y}) +\frac{\Delta}{4N^2}  }. \label{eqn:balnumdenom2}
\end{eqnarray}
In order to do this, we will appeal to the following identities:
\begin{eqnarray}
\dens_X+\dens_Y &=& \densca+ \denscu   \label{eqn:id1}  \\
\dens_{X,Y} &= &\frac{\densca+\denscu}{2}  \label{eqn:id2} \\
N \cdot  \densca + \Delta &=& N \cdot \denscu \label{eqn:id3}  \\
\densca & = & \dens_{X,Y}-\frac{\Delta}{2N} \label{eqn:id4}   \\
 \denscu & = & \dens_{X,Y}+\frac{\Delta}{2N}  \label{eqn:id5}
\end{eqnarray}
Equation (\ref{eqn:id1}) holds by simple inclusion-exclusion, Equation (\ref{eqn:id2}) follows directly from Equation (\ref{eqn:id1}),
Equation (\ref{eqn:id3}) is combinatorially trivial, and Equations (\ref{eqn:id4}) and (\ref{eqn:id5})
follow from Equations~(\ref{eqn:id2})~and~(\ref{eqn:id3}).

It is trivial to see that $\densca$ and $\dens_{X,Y}$ are balanced, so we need only compute
$\overline{\dens_X \dens_Y}$. Indeed, for any $h \in \stuff$ and and any
$(x,y) \in \x_h$, we have the following (using the identities
in Equations (\ref{eqn:id1}) through (\ref{eqn:id5}), and combinatorial symmetry, and
well-known identities involving binomial coefficients):
\begin{eqnarray}
\overline{\dens_X \dens_Y}(x,y) & = & \frac{1}{2^{\Delta (x,y)}} \sum_{(x',y') \in \x_h} \dens_X(x')\dens_Y(y') \nonumber \\
&=& \frac{1}{2^{\Delta (x,y)}} \sum_{i=0}^{\Delta (x,y)} {\Delta(x,y) \choose i} \ \ \frac{N \densca (x,y) +i }{N}
\ \ \frac{N \densca (x,y) + \Delta(x,y) -i }{N}\nonumber \\
&=& \frac{\densca (x,y) \denscu (x,y)}{2^{\Delta (x,y)}} \sum_{i=0}^{\Delta (x,y)} {\Delta(x,y) \choose i} \nonumber \\
& & \hspace{1in} + \ \ \frac{\Delta (x,y)}{2^{\Delta (x,y)}N^2}
\sum_{i=0}^{\Delta (x,y)} {\Delta(x,y) \choose i}i \ \ - \ \  \frac{1}{2^{\Delta (x,y)}N^2}
\sum_{i=0}^{\Delta (x,y)} {\Delta(x,y) \choose i}i^2 \nonumber  \\
& = & \densca (x,y) \denscu (x,y) \ \ + \ \ \frac{\Delta (x,y)}{2^{\Delta (x,y)}N^2} \left (  \Delta(x,y)2^{ \Delta(x,y)-1}  \right )\nonumber \\
& & \hspace{1in} - \ \ \frac{1}{2^{\Delta (x,y)}N^2} \left ( \Delta(x,y)+ \Delta^2(x,y) \right ) 2^{ \Delta(x,y)-2} \nonumber \\
& = & \left ( \dens_{X,Y}(x,y) - \frac{\Delta (x,y)}{2N}   \right ) \left ( \dens_{X,Y}(x,y) + \frac{\Delta (x,y)}{2N}   \right )
+ \frac{\Delta^2(x,y)}{2N^2} - \frac{\Delta(x,y)+\Delta^2(x,y)}{4N^2} \nonumber \\
& = & \dens_{X,Y}^2(x,y)-\frac{\Delta(x,y)}{4N^2}.  \label{eqn:punch}
\end{eqnarray}

Thus, by Equation (\ref{eqn:punch}) and the definition $\dens_{X,Y}=\frac{\dens_X+\dens_Y}{2}$ we obtain that
\begin{eqnarray} \overline{\densca - \dens_X \dens_Y}= \densca -\overline{\dens_X \dens_Y}=
\densca - \dens_{X,Y}^2 + \frac{\Delta}{4N^2}=\densca - \dens_X \dens_Y + \frac{1}{4} \left [ \frac{\Delta}{N^2} - (\dens_X-\dens_Y)^2 \right ] \mbox{and} \nonumber \\
\overline{\dens_{X,Y} - \dens_X \dens_Y}= \dens_{X,Y} - \overline{\dens_X \dens_Y}=
\dens_{X,Y} (1-\dens_{X,Y}) +\frac{\Delta}{4N^2} = \dens_{X,Y} - \dens_X \dens_Y + \frac{1}{4} \left [ \frac{\Delta}{N^2} - (\dens_X-\dens_Y)^2 \right ]; \ \ \  \label{eqn:forsd}
\end{eqnarray}
from this we have that Equation (\ref{eqn:balnumdenom2}), i.e. Equation (\ref{eqn:balnumdenom}), is proven, as desired. $\qed$ \\

\subsection{Proof of Corollary \ref{cor:cor} \label{sec:cor}}

Corollary \ref{cor:cor} states that
if the parameter space $\Theta$ is nondegenerate and also $\Theta \subseteq \paraoo$
then the statistic $\dens_{X,Y} \left ( 1- \dens_{X,Y} \right ) -  \frac{1}{2N} \left ( 1- \frac{1}{2N} \right ) \Delta $ is
UMVUE for $\sigma^2$. We prove this now.

Because $\Theta$ is nondegenerate, we pointed out in Example \ref{ex:bal} that
$\Delta: \x \rightarrow \R$ is the UMVU estimator of $2\sum_{i=1}^N(1-\varrho_i)p_i(1-p_i)$; here where
$\Theta \subseteq \paraoo$, we thus have that $\frac{\Delta}{2N}: \x \rightarrow \R$ is the UMVU estimator for
$\frac{1}{N}\sum_{i=1}^Np_i(1-p_i)=\mu(1-\mu)-\sigma^2$.

From Equation (\ref{eqn:forsd}) and Theorem \ref{thm:second} we have that
 $ \dens_{X,Y} (1-\dens_{X,Y})+ \frac{\Delta}{4N^2}$ is the UMVU estimator of
$\e \left ( \dens_{X,Y} - \dens_X \dens_Y \right )$,
which is equal to $\mu (1-\mu)$ since by hypothesis $\dens_X$ and $\dens_Y$ are independent.

Finally, by Corollary \ref{cor:vectorspace},
we have that  $ \dens_{X,Y} (1-\dens_{X,Y})+ \frac{\Delta}{4N^2} - \frac{\Delta}{2N}  $
is the UMVU estimator for $\mu (1-\mu) - [ \mu (1-\mu) -\sigma^2]$, and the result is shown. $\qed$

\subsection{Necessity of the $\Theta$ nondegeneracy assumption in Theorems \ref{thm:first}, \ref{thm:second} \label{sec:need} }

In the statement of Theorems \ref{thm:first} and \ref{thm:second}
we assume that the parameter space $\Theta$ is nondegenerate.
In this section we show that the claims of these theorems can fail if this condition is not satisfied.

Specifically, we will focus on a scenario in which the value of $\mu$ is known,
in which case the parameter space is reduced to
parameter tuples that have the prescribed value $\mu$.
This restricts the parameter space $\Theta$ to a particular hyperplane, which makes $\Theta$ degenerate; we will show
that the claims of Theorems \ref{thm:first} and \ref{thm:second} then fail, in general.

For simplicity, in this entire section, let us take $N=2$,
set $\varrho_1=\varrho_2=0$,
and suppose that $\mu:0<\mu<1$ is known; other than this, we allow $0<p_1<1$ and $0<p_2<1$.
Here, $\frac{p_1+p_2}{2}=\mu$ yields $p_2=2\mu-p_1$.
Denote $\delta:=\min \{ \mu, 1-\mu \}$, and denote $p:=p_1$; the parameter space is reduced to single
variable $p$ on the interval
$(\mu-\delta,\mu+\delta)$. There are $2^4=16$ points in $\x$; for each $(x,y) \in \x$, the probability of $(x,y)$ is given by
\[
\phi_{(x,y)}(p):=\left \{ \begin{array}{cl} p & \mbox{if } x_1=1 \\ 1-p & \mbox{if } x_1 =0 \end{array} \right \} \cdot
\left \{ \begin{array}{cl} p & \mbox{if } y_1=1 \\ 1-p & \mbox{if } y_1 =0 \end{array} \right \} \cdot
\left \{ \begin{array}{cl} 2\mu - p  & \mbox{if } x_2=1 \\ 1-2\mu+p & \mbox{if } x_2 =0 \end{array} \right \} \cdot
\left \{ \begin{array}{cl} 2\mu - p  & \mbox{if } y_2=1 \\ 1-2\mu+p & \mbox{if } y_2 =0 \end{array} \right \},
\]
which is a polynomial of degree $4$.

In Section \ref{sec:firstforward}, consider the linear system in Equation (\ref{eqn:specialcase});
that 9-by-9 linear system---describing statistic $S$ being an unbiased estimator of $g$---now becomes a 5-by-9 linear system over here.
This is because the columns of the left hand side matrix $A \otimes A$ and also the right hand side of the linear system, which there
were each 9-vectors consisting of the coefficients of particular polynomials in two variables,
can each now---in the reduced parameter space---be expressed as polynomials of degree $4$ in a single variable,
thus with five coefficients instead of $9$ coefficients.
As a 5-by-9 linear system, there is a nontrivial nullspace, and linear system solutions describing unbiasedness
are no longer unique, which implies that there will exist a statistic $T$ also an unbiased estimator of $g$ such that
it does not hold for all $h \in \stuff$ that
$\sum_{(x,y) \in \x_h} T(x,y) =  \sum_{(x,y) \in \x_h} S(x,y) $.
This completes our demonstration that the claim of  Theorem \ref{thm:first} may fail
in the absence of the nondegeneracy assumption for $\Theta$.\\

Next, we illustrate that the claim of Theorem \ref{thm:second} may fail without
the nondegeneracy assumption for $\Theta$.
The disagreement enumeration  statistic $\Delta$ is, by definition, unbiased for its expected
value $\e (\Delta)= 2\sum_{i=1}^N(1-\varrho_i)p_i(1-p_i)$, and we pointed out in Example \ref{ex:bal} that
$\Delta$ is the UMVU estimator of $\e (\Delta)= 2\sum_{i=1}^N(1-\varrho_i)p_i(1-p_i)$
if $\Theta$ is nondegenerate. Nonetheless, in general, we will see that there is no UMVU estimator
of $\e (\Delta)= 2\sum_{i=1}^N(1-\varrho_i)p_i(1-p_i)$ when $\mu$ is fixed.

%new
Indeed, consider $N=2$, any fixed value of $\mu$, and let the parameter space be 
parameterized by $p$ exactly as we did above in this section.
We will next formulate a quadratic program to find an unbiased estimator for
$\e (\Delta)$ which, for any given value of $p$, has least variance among the 
unbiased estimators. It will turn out that there is a unique solution to this
quadratic program. Then, if the solution differs for two different values of 
$p$ then there does not exists a UMVU estimator. Indeed, we performed computations, 
and found that this occurred when $\mu=.25$, as one example of many.

%previous
% in the
%setting described above. Indeed, for a specific example among many (for $N=2$, and parameter space parameterized by $p$ exactly
%as above in this section), take $\mu=.25$; we considered ALL unbiased estimators of $\e (\Delta)$ and then we considered
%which of them had the least variance for different fixed values of $p$.
%None of the unbiased estimators of $\e (\Delta)$ has variance less than
%or equal to the other unbiased estimators of $\e (\Delta)$, uniformly for all values of $p$.
%Thus there isn't a UMVU estimator for $\e (\Delta)$ here. From this we conclude that the claim of  Theorem~\ref{thm:second} may fail  without
%nondegeneracy of $\Theta$.\\

We now describe how to compute (in the scenario of this section) 
the unbiased estimator which, for any given $p$,  
has least variance among the unbiased estimators.

Let the points in $\x$ be ordered in any specified way, say $z_1,z_2,\ldots, z_{16}$.
Define the matrix $M \in \R^{5 \times 16}$ wherein, for all $i,j$, the entry $M_{ij}$ is the coefficient of $p^{i-1}$ in the
polynomial $\phi_{z_j}(p)$ (where $\phi_{z_j}(p)$ is as defined earlier
in this section). For any statistic $S:\x \rightarrow \R$, let $S$ be expressed as a vector $\vec{S} \in \R^{16}$
wherein, for all $i=1,2,\ldots,16$, we define $\vec{S}_i:=S(z_i)$. A function on the reduced parameter space
$g:(\mu-\delta,\mu+\delta) \rightarrow \R$ can only have an unbiased estimator if $g$ is a polynomial in the variable $p$
of degree at most $4$; this is  because $g$ would need to be
a linear combination of the $\phi$'s. Say that $\vec{g} \in \R^5$ is the vector wherein
for all $i=1,2,\ldots,5$, we define $\vec{g}_i$ to be the coefficient of $p^{i-1}$ in $g$.
Because $(\mu-\delta,\mu+\delta)$ is a nontrivial interval
(indeed, we just need at least 5 points) and by the uniqueness of interpolating polynomials,
we have that the unbiased estimators $S$ of any particular $g$ are precisely the solutions $\vec{S}$ of the linear system $M \vec{S}=\vec{g}$.

Suppose that there exists an unbiased estimator of $g$.
Among unbiased estimators of $g$, to find one of minimum variance {\bf for any specific value of} $p \in (\mu-\delta,\mu+\delta)$,
we proceed as follows. Let the vector of sample point probabilities for the respective $16$ sample space points
 be denoted $\vec{\varpi} \in \R^{16}$; we have
$\vec{\varpi}^T = [p^0,p^1,p^2,p^3,p^4]\cdot M$, which is a positive vector since $p \in (\mu-\delta,\mu+\delta)$;
finding a (globally) unbiased estimator with (specifically for $p$) minimum variance is equivalent to the
quadratic, convex optimization problem
$\min \sum_{i=1}^{16} \vec{\varpi}_i\vec{S}_i^2$ such that $\vec{S}$ satisfies
$M \vec{S}=\vec{g}$ (we minimize the second moment for the estimator, since the first moment is fixed).
Define a bijective change of variables where new variables $\vec{S}' \in \R^{16}$ are such that for all $i=1,2,\ldots,16$ we
have $\vec{S}'_i:=\sqrt{\vec{\varpi}}_i \vec{S}_i$, and define $M' \in \R^{5 \times 16}$ such that,
for all $i,j$, we have $M'_{ij}=\frac{1}{\sqrt{\vec{\varpi}_j}}M_{ij}$. Now this minimum variance problem is equivalent to
$\min \| \vec{S}' \|_2$ such that $\vec{S}'$ satisfies
$M' \vec{S}'=\vec{g}$. Classical generalized inverse theory guarantees a unique solution $\vec{S}'=M'^\dag \vec{g}$ (the symbol $\dag$
denotes the Moore-Penrose Generalized Inverse of the matrix), which corresponds to
statistic $S$ wherein $S(z_i)=\frac{(M'^\dag \vec{g})_i}{\sqrt{\vec{\varpi_i}}}$ for each $i=1,2,\ldots,16$, which is unique as having minimum variance
(for the particular value of $p$) among the (globally) unbiased estimators.

This concludes the description of the way we computed, in the scenario of this section, the unbiased estimator which, for a fixed value of $p$,
has least variance. (An excellent reference for matrix analysis in general, with theory
of generalized inverses, is \cite{HJ}.)

\section{Simulation experiments: comparing $\st$,  $\overline{\st}$, and $\st'$ \label{sec:empir}}

As we mentioned earlier, in \cite{TotCor} it was empirically demonstrated---for correlated graphs in broad families within
our model---that graph matching complexity as well as graph matchability are each functions of total correlation,
and it was also proved in \cite{TotCor} that alignment strength $\st$ is a strongly consistent estimator of $\rho_T$.
The specific formulation/definition of alignment strength $\st$ arose in a very natural way; see \cite{TotCor}.
Nonetheless, $\st$ suffers from a deficiency; in Example \ref{ex:nonobal} we pointed out that $\st$ is not balanced.
The balanced statistic $\overline{\st}$ reduces the variance, keeping the expected value unchanged.
In this section we will empirically demonstrate that another balanced statistic, denoted $\st'$, is often superior to
$\overline{\st}$ in estimating $\rho_T$. Note that there is no contradiction to Theorem \ref{thm:second}, which
asserts that, assuming the parameter space is nondegenerate, $\overline{\st}$ is UMVUE for $\e (\st)$; indeed,
$\st'$ can be biased as an estimator of $\e (\st)$. Which can be a good thing; we will see that $\st'$ frequently
has less bias than $\overline{\st}$ in the estimation of $\rho_T$, and in all of these experiments here
$\st'$ has less variance than $\overline{\st}$.

But we first make a computationally helpful observation about computing the value of
$\overline{\st}$.\\

In general, when given an arbitrary  statistic $S: \x \rightarrow \R$, the computation of the value of $\overline{S}$, even for just one
particular sample space point $(x,y) \in \x$, can require exponential time; indeed,
there are $2^{\Delta(x,y)}$ values to average. In the case of computing $\overline{\st}$, this computation
can be greatly simplified as follows. Given any particular $h \in \stuff$ and any
particular $(x,y) \in \x_h$ (such that not both $x$ and $y$ are all zeros, and not both $x$ and $y$ are all ones), we have
by Equation (\ref{def:asother}) and Equation (\ref{eqn:id3}), that
\begin{eqnarray*}
\overline{\st}(x,y) & := &  \frac{1}{2^{\Delta(x,y)}}\sum_{(x',y')\in \x_h}\st(x',y')   \nonumber \\
&=& \frac{1}{2^{\Delta (x,y)}} \sum_{i=0}^{\Delta (x,y)} {\Delta(x,y) \choose i}
\frac{\densca (x,y) - \Big (\densca(x,y) + \frac{i}{N} \Big ) \left ( \densca(x,y)+\frac{\Delta(x,y)}{N}-\frac{i}{N} \right ) }{
\dens_{X,Y}(x,y) -\Big (\densca(x,y) + \frac{i}{N} \Big ) \left (\densca(x,y)+\frac{\Delta(x,y)}{N}-\frac{i}{N} \right ) }
\nonumber \\
\end{eqnarray*}
The above provides a linear time
algorithm for computing $\overline{\st}$ for any $(x,y) \in \x$, although this computation is much more involved then the
very simple formula for $\st'$ as given in Remark \ref{rem:quasi}.\\

Now we are prepared to do a simulation experiment to compare the variances of $\st$, and $\overline{\st}$, and $\st'$.
The expected values of $\st$ and $\overline{\st}$ are of course the same, so there is no difference between their biases
in estimating $\varrho_T$. However, the expected values of $\st$ and $\st'$ are not the same, so we want to also
compare their biases in estimating $\varrho_T$.

In the first set of experiments, we did $200$ independent replicates of the following experiment.
We realized $p_1,p_2,p_3,p_4,p_5,p_6,\varrho_1,\varrho_2,\varrho_3,\varrho_4,\varrho_5,\varrho_6$
(which correspond to the six pairs of vertices in a vertex set with four vertices)
independently from a Uniform$(0,1)$ distribution; the first five such experiments' values were:
\[
\begin{array}{cccccccccccc}
       p_1 &    p_2   &    p_3   &    p_4   &     p_5  &     p_6  &  \varrho_1 &  \varrho_2 &  \varrho_3 &  \varrho_4 &  \varrho_5 &  \varrho_6 \\ \hline
    0.6892 &   0.7224 &   0.4795 &   0.8985 &   0.4022 &   0.7043 &   0.8429 &   0.9852  &  0.8006 &   0.3118 &   0.5768 &   0.5751 \\
    0.7482 &   0.1499 &   0.6393 &   0.1182 &   0.6207 &   0.7295 &   0.8988 &   0.6088  &  0.7388 &   0.0553 &   0.9440 &   0.0100 \\
    0.4505 &   0.6596 &   0.5447 &   0.9884 &   0.1544 &   0.2243 &   0.9390 &   0.2537  &  0.1417 &   0.7538 &   0.8715 &   0.8094 \\
    0.0838 &   0.5186 &   0.6473 &   0.5400 &   0.3813 &   0.2691 &   0.8154 &   0.1326  &  0.4379 &   0.1319 &   0.5076 &   0.6088 \\
    0.2290 &   0.9730 &   0.5439 &   0.7069 &   0.1611 &   0.6730 &   0.0014 &   0.5450  &  0.3504 &   0.3559 &   0.7888 &   0.4799
\end{array}
\]
Then, for these realized parameters, we computed
(exactly, by enumerating the sample space and sample point probabilities) the values of
$\e (\st)$, $\e (\st')$, $\varrho_T$, Var$(\st)$, Var$(\overline{\st})$, and Var$(\st')$. Of course,
$\e(\st)=\e(\overline{\st})$. The first five experiments' outcomes were:
\[
\begin{array}{ccc|ccc}
      \e(\st) & \e (\st') & \varrho_T & \textup{Var}(\st) & \textup{Var} (\overline{\st}) & \textup{Var}(\st') \\ \hline
    0.6851 &   0.6857 &   0.7516 &   0.1219 &   0.1214  &  0.1206 \\
    0.6835 &   0.6843 &   0.7093 &   0.0885 &   0.0879  &  0.0870 \\
    0.6827 &   0.6833 &   0.7011 &   0.0745 &   0.0740  &  0.0734 \\
    0.4310 &   0.4339 &   0.4697 &   0.1345 &   0.1318  &  0.1291 \\
    0.5619 &   0.5635 &   0.5789 &   0.1073 &   0.1059  &  0.1043
\end{array}
\]
In every one of these $200$ experiments, we had that
$ \textup{Var}(\st) > \textup{Var} (\overline{\st}) > \textup{Var}(\st')$. In $199$ of these $200$
experiments we had  $\e(\st) < \e (\st') < \varrho_T$.

We then repeated the $200$ experiments, except that $\varrho_i=0$ for all $i=1,2,3,4,5,6$.
The first five experiments' outcomes were:
\[
\begin{array}{ccc|ccc}
      \e(\st) & \e (\st') & \varrho_T & \textup{Var}(\st) & \textup{Var} (\overline{\st}) & \textup{Var}(\st') \\ \hline
    0.3278 &   0.3320 &   0.3234 &   0.1222 &   0.1182 &   0.1149 \\
    0.4965 &   0.4986 &   0.4918 &   0.1052 &   0.1033 &   0.1014 \\
    0.4240 &   0.4269 &   0.4169 &   0.1098 &   0.1072 &   0.1048 \\
    0.1260 &   0.1335 &   0.1333 &   0.1433 &   0.1354 &   0.1307 \\
    0.5204 &   0.5225 &   0.5177 &   0.1094 &   0.1076 &   0.1056 \\
\end{array}
\]
Again, in every one of these $200$ experiments, we had that
$ \textup{Var}(\st) > \textup{Var} (\overline{\st}) > \textup{Var}(\st')$.
However, in only $41$ of the $200$ experiments was the bias of $\st'$
less than that of $\st$, meaning that $|\e(\st')-\varrho_T|< |\e(\st)-\varrho_T|$.

We then repeated the first $200$ experiments, except that $p_i=\frac{1}{2}$ for all $i=1,2,3,4,5,6$.
The first five experiments' outcomes were:
\[
\begin{array}{ccc|ccc}
      \e(\st) & \e (\st') & \varrho_T & \textup{Var}(\st) & \textup{Var} (\overline{\st}) & \textup{Var}(\st') \\ \hline
    0.6866 &   0.6872  &  0.7392 &   0.0989  &  0.0983  &  0.0976  \\
    0.3062 &   0.3108  &  0.3496 &   0.1375  &  0.1330  &  0.1293  \\
    0.3776 &   0.3812  &  0.4257 &   0.1367  &  0.1333  &  0.1301  \\
    0.3745 &   0.3781  &  0.4221 &   0.1384  &  0.1349  &  0.1316 \\
    0.6384 &   0.6393  &  0.6919 &   0.1095  &  0.1086  &  0.1075
\end{array}
\]
Again, in every one of these $200$ experiments, we had that
$ \textup{Var}(\st) > \textup{Var} (\overline{\st}) > \textup{Var}(\st')$.
In all these $200$
experiments we had  $\e(\st) < \e (\st') < \varrho_T$. (In all of the above
$600$ experiments, we adopted the convention that the statistics have the value $0$
at the two sample points $(x,y)$ where $x$ and $y$ are all zeros and where they are all ones.
Indeed, we saw empirically that this choice had
a negligible numerical impact on the experiments here.)

Remarkably, in all of the many tens of thousands of experiments that we conducted, for many different parameter values,
we always found that the mean squared
error in estimating $\rT$, denote it $MSE( \cdot, \rT)$, was lower for the modified alignment strength $\st'$ than for
the balanced alignment strength $\overline{\st}$. Based on these computations, we conjecture the following.

\begin{conjecture} For all $N$ and $\theta \in \para$, it holds that $MSE(\st',\rT) \leq MSE(\overline{\st},\rT)$. \label{conje}
\end{conjecture}

\section{Summary and future directions}

Our setting is the correlated Bernoulli random graph model
 for the production of a pair of correlated random graphs,
wherein different pairs of vertices are allowed different probabilities
of adjacency, and inter-graph edge correlations are allowed to be
different for different pairs of vertices. This is a broad and useful model.
Our main results come in two groups.

The first group of results: We introduce a “balancing” procedure to lower
the mean squared error for any statistic used to estimate any function of
the model parameters; it is essentially a Rao-Blackwellization procedure
utilizing the disagreement vector statistic, which we prove is complete and sufficient.
Indeed, given any unbiased estimator of any function of the model parameters,
we neatly characterize all unbiased estimators, as well as the UMVUE estimator for
this function of the model parameters. With these tools, we obtain the
second group of results, which involve estimating the total correlation
parameter, which is of current interest in the theory of Graph Matching,
and has been recently shown to play a critical role in matchability and also in  graph matching
runtime complexity (when graph matching is solved exactly via integer programming).\cite{TotCor}

Future steps would be to extend our results in this paper to
broader random graph models and to settle Conjecture \ref{conje}.\\ \\

\noindent {\bf Acknowledgements} This material is based on research sponsored by the Air Force Research Laboratory and Defense Advanced Research Projects Agency (DARPA) under agreement number FA8750-20-2-1001 and FA8750-18-2-0035.
This work is also supported in part by the D3M program of DARPA.
The U.S. Government is authorized
to reproduce and distribute reprints for Governmental purposes notwithstanding any
copyright notation thereon. The views and conclusions contained herein are those of
the authors and should not be interpreted as necessarily representing the official policies
or endorsements, either expressed or implied, of the Air Force Research Laboratory and
DARPA or the U.S. Government. The authors also gratefully acknowledge the support
of NIH grant BRAIN U01-NS108637. The authors also
thank an anonymous referee for very helpful and
thoughtful comments that added much to the quality of
the exposition.


\begin{thebibliography}{9}

\bibitem{Blackwell}
    D.~Blackwell,
    Conditional expectation and unbiased sequential estimation,
    {\em The Annals of Mathematical Statistics} {\bf 18}, (1947), pp 105-–110.

\bibitem{thirtyyears}
    D.~Conte, P.~Foggia, C.~Sansone, M.~Vento,
    Thirty years of graph matching in pattern recognition,
    {\em Int. J. Pattern Recognit. Artif. Intell.} {\bf 18} (2004), pp 265-–298.

\bibitem{SGM1}
    D.E.~Fishkind, S.~Adali, H.G.~Patsolic, L.~Meng, D.~Singh, V.~Lyzinski, C.E.~Priebe,
    Seeded Graph Matching,
    {\em Pattern Recognition}, {\bf 87} (2019), pp 203--215.

\bibitem{TotCor}
    D.E.~Fishkind, L.~Meng, A.~Sun, C.E.~Priebe, V.~Lyzinski,
    Alignment strength and correlation for graphs,
    {\em Pattern Recognition Letters}, {\bf 125} (2019), pp 295--302.

\bibitem{tenyears}
    P.~Foggia, G.~Perncannella, M.~Vento,
    Graph matching and learning in pattern recognition in the last 10 years,
    {\em Int. J. Pattern Recognit. Artif. Intell.} {\bf 28} (2014).

\bibitem{HJ}
    R.A.~Horn and C.R.~Johnson,
    {\em Matrix Analysis}, Cambridge University Press, (1985), 2nd edition (2013).

\bibitem{Lang}
    S.~Lang, {\em Algebra}, Springer, 3rd Edition (2005).

\bibitem{LSone}
    E.L.~Lehmann and H.~Scheffe,
    Completeness, similar regions and unbiased estimation:~Part~I,
    {\em Sankhya: The Indian Journal of Statistics} {\bf 10}, (1950), pp 305–-340.

\bibitem{LStwo}
    E.L.~Lehmann and H.~Scheffe,
    Completeness, similar regions and unbiased estimation:~Part~II,
    {\em Sankhya: The Indian Journal of Statistics} {\bf 15}, (1955), pp 219–-236.

\bibitem{GMrelax}
    V.~Lyzinski, D.E.~Fishkind, M.~Fiori, J.T.~Vogelstein, C.E.~Priebe, G.~Sapiro,
    Graph matching: Relax at your own risk,
    {\em IEEE Transactions on Pattern Analysis and Machine Intelligence}, {\bf 38} (2016), pp 60--73.

\bibitem{SGM2}
    V.~Lyzinski, D.E.~Fishkind, C.E.~Priebe,
    Seeded graph matching for correlated Erdos-Renyi graphs,
    {\em Journal of Machine Learning Research} {\bf 15} (2014), pp 3513--3540.

\bibitem{Rao}
    C.R.~Rao,
    Information and accuracy attainable in the estimation of statistical parameters,
    {\em Bulletin of the Calcutta Mathematical Society} {\bf 37}, (1945), pp 81-–91.

\bibitem{gmsurvey}
    J.~Yan, X.C.~Yin, W.~Lin, C.~Deng, H.~Zha, X.~Yang,
    A short survey of recent advances in graph matching,
    {\em Proceedings of the 2016 ACM on International Conference
    on Multimedia Retrieval, ACM} (2016), pp 167-–174.

\end{thebibliography}
\end{document}